\documentclass[twoside,reqno,final, 11pt]{amsart}
\usepackage{amsmath}
\usepackage{amssymb}
\usepackage{graphicx}
\usepackage{mathrsfs}
\usepackage{amsrefs}

\input epsf

\numberwithin{equation}{section}
\newtheorem{Thm}{Theorem}
\newtheorem{Cor}{Corollary}

\newtheorem{Lem}{Lemma}
\newtheorem{Conj}{Conjecture}
\theoremstyle{definition}
\newtheorem{Def}{Definition}

\theoremstyle{remark}
\newtheorem{Rem}{Remark}

\renewcommand\ge{\geqslant}
\renewcommand\le{\leqslant}
\let\tildeaccent=\~
\let\hataccent=\^
\renewcommand\~[1]{\widetilde{#1}}
\renewcommand\^[1]{\widehat{#1}}

\def\<{\left<}
\def\>{\right>}
\def\({\ifmmode\left(\else\textup{(}\fi}
\def\){\ifmmode\right)\else\textup{)}\fi}

\def\e{\varepsilon}
\def\R{{\mathbb R}}
\def\C{{\mathbb C}}
\def\N{{\mathbb N}}
\def\Z{{\mathbb Z}}
\def\l{\lambda}
\def\Re{\operatorname{Re}}
\def\Im{\operatorname{Im}}
\def\const{\operatorname{const}}
\def\M{{\mathscr M}}
\def\L{\varLambda}

\let\epsilon=\varepsilon
\let\kappa=\varkappa

\def\setD{\mathscr{D}}
\def\setDC{\mathscr{D}_\C}
\def\setB{\mathscr{B}}
\def\setBC{\mathscr{B}_\C}

\def\spaceJ{\mathfrak{J}}
\def\spaceO{\mathfrak{O}}

\def\operC{\mathcal{C}}
\def\operE{\mathcal{E}}
\def\operD{\mathcal{D}}
\def\operP{\mathcal{P}}
\def\operK{\mathcal{K}}

\begin{document}

\title [Perturbations of Darbouxian systems]
 { On limit cycles appearing by polynomial perturbation of
                  Darbouxian integrable systems }

\author[D.~Novikov]
 {Dmitry Novikov$^*$ }

\email {\tt dnovikov@wisdom.weizmann.ac.il}

\date{\today}


\begin{abstract}
 We prove an existential finiteness result for integrals of
 rational 1-forms over the level curves of Darbouxian integrals.
\end{abstract}

 \maketitle

\section{Limit cycles born by perturbation of integrable systems}

\subsection{Poincar\'e--Pontryagin integral}
Limit cycles (isolated periodic trajectories) of polynomial planar
vector fields can be produced by perturbing integrable
systems which have nested continuous families of non-isolated
periodic trajectories. The number and position of limit cycles
born in such perturbations is determined by the number and
position of \emph{isolated zeros} of the Poincar\'e--Pontryagin
integral of the dissipation form along the closed periodic orbits
of the non-perturbed integrable vector field. The
\emph{infinitesimal Hilbert problem} is to place an upper bound
for the number of isolated zeros of this integral. This problem
was repeatedly formulated as a relaxed form of the Hilbert 16th
problem by V.~Arnold \cite{arnold:problems}.

Instead of polynomial vector fields, it is more convenient to deal
with (singular) foliations of the real plane $\R^2$ by solutions
of \emph{rational} Pfaffian equations $\theta=0$, where $\theta$
is a 1-form on $\R^2$ with rational coefficients (since only the
distribution of null spaces of the form $\theta$ makes geometric
sense, one can always replace the rational form by a polynomial
one). The foliation $\mathcal F$ is (Darbouxian)
\emph{integrable}, if the form is \emph{closed}, $d\theta=0$. An
integrable foliation always has a ``multivalued'' first integral
of the form $f(x,y)=\exp r(x,y)\cdot\prod_j p_j(x,y)^{\l_j}$,
where $r(x,y)$ is a rational function and $p_j$ polynomials in
$x,y$, which are involved in the in general non-rational powers
$\l_j$. A particular case of the integrable foliations consists of
\emph{Hamiltonian foliations} defined by the \emph{exact
polynomial} 1-form $\theta=dh$, $h\in\R[x,y]$.

If $L\subset\R^2$ is a compact smooth leaf (oval) of an integrable
foliation $\theta=0$ belonging to the level curve
$\{f=a\}\subset\R^2$, then all nearby leaves of this foliation are
also closed and form a continuous family of ovals belonging to the
level curves $\{f=t\}$, $t\in(\R^1,a)$. The ovals from this
family, denoted by $L_t$, are uniquely parameterized by the values
of the real variable $t$ varying in an open interval
$\alpha<t<\beta$. Their union is an annulus $A_{\alpha,\beta}$
bounded by the Hausdorff limits $\overline L_\alpha=\lim_{t\to
\alpha+0} L_t$ and $\overline L_\beta=\lim_{t\to\beta-0}$. The
limits themselves may be not leaves of the foliation, but rather
\emph{separatrix polygons} (unions of several leaves and one or
more singular points of the foliation $\mathcal F$). The polar
locus of the form $\theta$ is usually a separatrix polygon.

Consider an one-parameter family of planar real foliations
$\mathcal F_\e$ defined by the Pfaffian equations
\begin{equation}\label{pe}
    \theta+\e\omega=0,\qquad
    \operatorname{Poles}(\omega)\subseteq\operatorname{Poles}(\theta),
\end{equation}
with a closed rational 1-form $\theta$ and arbitrary rational 1-form
$\omega$. The assumption on the poles guarantees that the
perturbation \eqref{pe} does not create additional singularities in
the annulus $A_{\alpha,\beta}$ (such $\omega$ will be called
\emph{admissible} for $\theta$). The \emph{Poincar\'e--Pontryagin
integral}\footnote{The function $I(t)$ is sometimes called Melnikov
function or (misleadingly) the Abelian integral, see below.}
associated with this family, is the function
\begin{equation}\label{ai}
    t\mapsto I(t)=\oint_{L_t}\omega,\qquad
    L_t\subseteq\{f=t\}\text{ an oval, }
    t\in (\alpha,\beta).
\end{equation}
The following well-known \emph{Poincar\'e--Pontryagin criterion}
establishes connections between the integral \eqref{ai} and limit
cycles of the system \eqref{pe}: \emph{If the function $I(t)$ has
an isolated zero of multiplicity $\mu$ at an interior point
$a\in(\alpha,\beta)$, then the perturbed foliation $\mathcal F_\e$
for all sufficiently small values of $\e$ and $\delta>0$ has no
more than $\mu$ limit cycles in the thin annulus between the ovals
$L_{a-\delta}$ and $L_{a+\delta}$}. The assertion is reasonably
sharp: under certain additional assumptions of nondegeneracy one
can produce exactly $\mu$ limit cycles.

Thus bounds for the number of isolated zeros for the integral
\eqref{ai} translate into bounds for the number of limit cycles
for near-integrable systems.\footnote{However, the above
Poincar\'e--Pontryagin criterion does not apply to limit cycles
born from separatrix polygons $\overline L_\alpha,\overline
L_\beta$ which bound the annuli of periodic orbits of the
integrable systems. This bifurcation requires much deeper analysis
that will not be discussed here.} The following conjecture seems
to be believed by most experts in the area.

\begin{Conj}\label{conj:gen}
For any pair of rational 1-forms $(\theta,\omega)$ on the real
plane $\R^2$ of degrees $n,m$ respectively, such that $d\theta=0$
and
$\operatorname{Poles}(\omega)\subseteq\operatorname{Poles}(\theta)$,
the number of isolated real zeros of the integral \eqref{ai} is
bounded by a constant $N=N(n,m)$ depending only on $n$ and $m$.
\end{Conj}

Demonstration (or refutation) of this Conjecture constitutes the
\emph{infinitesimal Hilbert Sixteenth problem}
\cite{arnold:problems,centennial}. A separate question concerns
computability (and practical computation) of the bound $N(n,m)$.

\subsection{Analytic properties of the Poincar\'e--Pontryagin integral:
Ham\-iltonian versus Darbouxian cases} The function $I(t)$ defined
by the integral \eqref{ai} is obviously real analytic inside the
real interval $(\alpha,\beta)$, and hence the number of isolated
zeros is necessarily finite on any compact subinterval of this
interval. Yet this circumstance neither implies any explicit upper
bound, nor rules out the accumulation of infinitely many isolated
roots to the boundary points $\alpha,\beta$ of the interval.

Most studies of the function $I(t)$ are based on its analytic
continuation from the real interval $(\alpha,\beta)$ to the
complex domain. Already at this step drastic differences between
Hamiltonian and general integrable cases occur.

\subsubsection{Abelian integrals}
If the initial foliation $\mathcal F=\{\theta=0\}$ is Hamiltonian,
i.e., if $\theta=dh$, $h\in\R[x,y]$, then $\omega$ must also be
polynomial and the function $I(t)$ is an \emph{Abelian integral},
integral of a polynomial 1-form over a cycle on \emph{algebraic
curve}. Such integrals possess the following properties.

(i) The integral $I(t)$ admits analytic continuation on the entire
 Riemann sphere
 $\C P^1=\C\cup\{\infty\}$ as a multivalued function ramified over
 finitely many points. For a generic Hamiltonian $h$ these points
 are exactly the (complex) critical values of $h$.

(ii) $I(t)$ is (a linear combination of coordinates of) a
 solution of a Fuchsian linear system on $\C P^1$.

(iii) Near each singular point $t_0$ of the Fuchsian system,
$I(t)$
 admits a local representation as the finite sum $I(t)=\sum _{j,k}
 a_{jk}(t)\,(t-t_0)^{\l_j}\ln^k(t-t_0)$ with the coefficients $a_{jk}$
 holomorphic at the point $t_0$ and finitely many rational exponents
 $\l_j\in\mathbb Q$.

The last property alone suffices to guarantee that Abelian integrals
are \emph{non-oscillating}: any such integral has only finitely many
real isolated zeros on any interval. This can be proved by
application of the ``derivation-division algorithm''
\cite{roussarie-book}. Moreover, the Fewnomial theory
\cite{fewnomials} which can be regarded as a multidimensional
generalization of the above algorithm, implies the following uniform
version of the non-oscillation property, achieved in
\cite{varchenko,khovanskii}.

\begin{Thm}[A.~Khovanskii--A.~Varchenko, 1984]\label{thm:kv}
The total number of real isolated zeros of Abelian integrals over
the level curves of a polynomial $h$ is uniformly bounded over all
polynomial 1-forms $\omega$ and all Hamiltonians $h$ sufficiently
close to any given combination $(h_0,\omega_0)$.
\end{Thm}

By the standard compactness and homogeneity arguments,
Theorem~\ref{thm:kv} implies that for any combination of
$n,m\in\N$ the number of real isolated zeros of Abelian integrals
of forms of degree $\le m$ over Hamiltonians of degree $\le n$ is
bounded by a constant $N_0(n,m)$ depending only on $n,m$. However,
this constant is absolutely existential (non-constructive). Other,
completely different tools allow to place  \emph{explicit} upper
bounds on the number of zeros of Abelian integrals which are at
least $\delta$-distant from the critical values of the Hamiltonian
\cite{meandering,quasialg} (the bounds depend on $\delta>0$).

One can expect similar behavior of the Poincar\'e--Pontryagin
integrals in the case when $\theta$ is exact albeit non-polynomial
anymore, $\theta=df$, with $f(x,y)=\frac{P(x,y)}{Q(x,y)}$ a rational
function, though the subject was not explored systematically.

\subsubsection{Pseudoabelian integrals}
Completely different is the picture if the (unperturbed) closed form
$\theta$ is \emph{non-exact}. In this article we consider only the
simplest case when $\theta$ has at most first order poles after
extension on $\C P^2$ and hence a Darbouxian integral of the form
\begin{equation}\label{dai}
    f(x,y)=\prod_j p_j(x,y)^{\l_j},\qquad \l_j\in\R_+,\
    j=1,\dots,m,
\end{equation}
ramified over the separatrix polygon
$S(\theta)=\bigcup_j\{p_j=0\}\subset\R^2$ corresponding to the
``critical'' value $f=0$. For lack of better name, we will refer
to integrals of the form
\begin{equation}\label{eq:pa}
    I(t)=\oint_{L_t}\omega,\qquad
    L_t\subseteq\{f=t\}\subseteq\R^2,
\end{equation}
where $\omega$ is a rational form without singularities on $L_t$,
as \emph{pseudoabelian integrals}.

This particular case already contains all difficulties which
distinguish the general integrable case from the Hamiltonian one.
The closed level curves $\{f=t\}$ cease to be algebraic unless all
ratios $\l_j/\l_k$ are rational, moreover, their complexifications
are generically dense in $\C^2$ \cite{thebook}. The topological
methods of continuation of the corresponding integral to the complex
domain, which work so nicely in the Hamiltonian case, fail
\cite{paul}, thus the analytical continuation of the corresponding
integral $I(t)$ should be achieved by completely different methods.
Besides, $I(t)$ is not known to satisfy any reasonable linear or
nonlinear differential equation of finite order, which makes
impossible applications of the methods from
\cite{meandering,quasialg}. For example, the local representation
\eqref{moura} already implies that $I(t)$ is not a solution of a
Fuchsian ordinary differential equation of finite order.

Even the most basic among the properties of the Abelian integrals,
the local analytic representation, fails for pseudoabelian
integrals. More precisely, it can be shown that if the separatrix
curve $S(\theta)$ has only normal crossings and the collection of
Darboux exponents $[\l_1:\cdots:\l_n]$ satisfies certain {generic}
\emph{arithmetic} properties, then the pseudoabelian integral
$I(t)$ can be locally near $t=0$ represented as a composition,
\begin{equation}\label{moura}
    I(t)=F(t,t^{\mu_1},\cdots,t^{\mu_n}), \qquad
    \mu_1,\dots,\mu_n\in\R_+,
\end{equation}
with a real analytic function $F$ of $n+1$ variables. This
observation, due to C.~Moura \cite{moura}, suffices to prove
non-oscillatory behavior of pseudoabelian integrals, but under the
above arithmetical restrictions on the Darboux exponents. Such form
of arithmetics-conditioned finiteness is not unknown in the general
theory of $o$-minimal structures, see \cite{soufflet}. Clearly, if
this arithmetic restriction is indeed necessary to rule out
non-oscillation of an individual pseudoabelian integral $I(t)$,
there would be almost no hope for uniform bounds for the number of
zeros of pseudoabelian integrals.

\subsection{Principal result}
Our main result is the \emph{uniform non-oscillation} of
pseudoabelian integrals associated with a \emph{generic} Darbouxian
integrable foliation near the separatrix polygon $S(\theta)$. To
formulate it accurately, we introduce Darbouxian classes of
integrable systems. Let $m_1,\dots,m_n$ be a tuple of natural
numbers, $m_j\in\N$.

\begin{Def}
The Darbouxian class $\mathscr D=\mathscr D(m_1,\dots,m_n)$ is the
class of real integrable foliations on $\R P^2$ defined by the
Pfaffian equations $\theta=0$, where $\theta$ is a
\emph{logarithmic form} with constraints on the degrees as
follows\footnote{One can easily show (e.g., see \cite{thebook})
that a closed rational form $\theta$ with poles of order $\le 1$
on $\C P^2$ has the form \eqref{log-form} with complex constants
$\l_j$.},
\begin{equation}\label{log-form}
    \theta=\sum_{j=1}^n\l_j\frac{dp_j}{p_j},
    \qquad p_j\in\R[x,y],\ \deg p_j\le m_j, \ \l_j>0.
\end{equation}
The polynomials $p_j$ are  assumed to be irreducible and different,
hence the polar locus $S(\theta)$ is the union of $n$ algebraic
curves $\{p_j=0\}\subset\R^2$.
\end{Def}

Foliations from a fixed Darbouxian class are parameterized by
points from an open subspace in the Euclidean space with the
coordinates being the exponents $\l_1,\dots,\l_n>0$ and the
respective polynomials $p_j$, identified with their coefficients.
Using this identification, one can define open neighborhoods of a
given closed form $\theta=0$ in its Darbouxian class $\mathscr D$.
We will consider only foliations which have continuous families of
ovals accumulating to the separatrix $S(\theta)$ or its proper
subset (as usual, in the Hausdorff sense).

If a real rational 1-form $\omega$ on $\R^2$ is \emph{admissible}
for the logarithmic form $\theta$, i.e., if the polar locus of
$\omega$ belongs to the separatrix $S(\theta)$, then the integral
$I(t)=\oint_{L_t}$ over any non-singular real oval $L_t$, $t\ne
0$, is well defined. We claim that the number of real isolated
zeros of such integrals is locally uniformly bounded.

\begin{Thm}[uniform non-oscillation of pseudoabelian
 integrals]\label{thm:main}
Let the form $\theta_0=\sum_{j=1}^n\l_j\frac{dp_j}{p_j}$ defining the
Darbouxian integrable foliation $\mathcal F=\{\theta_0=0\}$ be from
the class $\mathscr D(m_1,\dots,m_n)$, and assume that the real
algebraic curves $\{p_j=0\}$ are smooth and intersect transversally
(in particular, there are no triple intersections). Let $\omega_0$ be
a rational one-form admissible for $\theta_0$ with poles of order at
most $m_0$.

Then the number of isolated zeroes of the pseudoabelian integral near
$t=0$, as a function of $(\theta, \omega)$, is locally bounded at
$(\theta_0,\omega_0)$ over all admissible pairs $(\theta,\omega)$.
\end{Thm}

In other words, for any admissible pair of rational 1-forms
$(\theta_0,\omega_0)$ as in the Theorem, there exist a finite number
$N=N(\theta_0,\omega_0)$ and $\e>0$, depending on the pair, such that
the number of real isolated zeros of any pseudoabelian integral
corresponding to $(\theta,\omega)$ and to a bounded family of cycles
$L_t$, $t\in (0,\e)$, is no greater than $N$ provided that the pair
$(\theta,\omega)$ is sufficiently close to $(\theta_0,\omega_0)$.

After this result one may expect that the number of zeros of
pseudoabelian integrals would be uniformly bounded over all
Darbouxian classes. The following Conjecture is a relaxed form of
the Infinitesimal Hilbert 16th Problem (Conjecture~\ref{conj:gen})

\begin{Conj}
For any Darbouxian class $\mathscr D(m_1,\dots,m_n)$ and any
degree $m_0$ there exists a finite number $N=N(m_0,m_1,\dots,m_n)$
such that any pseudoabelian integral of a rational 1-form $\omega$
of degree $\le m_0$ along ovals of the corresponding Darbouxian
integral \eqref{dai} has no more than $N$ real isolated zeros.
\end{Conj}

Again in contrast with the Hamiltonian case, the Conjecture does
not follow from Theorem~\ref{thm:main}, since the compactness
arguments fail for pseudoabelian integrals. Indeed, the assumption
on the separatrix is an open condition in $\mathscr D$ which in
the limit may degenerate into non-transversal singularities at the
corners. Besides, the Darbouxian classes themselves are
non-compact and in the limit contain closed forms with poles of
higher orders and ``vanishing separatrices'' as $\l_j\to0^+$. All
these scenarios require additional study.

\subsection{Proofs: Mellin transform and generalized Petrov operators}
The most fundamental tool used to prove nonoscillation-type results,
is the classical derivation--division algorithm \cite{roussarie-book}
which goes as far back as to Descartes. In the modern language it
consists in constructing a differential operator $\operD$ (in
general, with variable coefficients) such that the result of the
derivation $\operD I(t)$ is a function without real isolated zeros
(e.g., identical zero). By the Rolle-Descartes theorem, the number of
zeros of zeros of $I$ is then bounded by a constant explicitly
expressible in terms of $\operD$ (in the simplest case of operators
with constant coefficients, by the order of $\operD$, provided that
all characteristic roots are real). This observation was generalized
in several directions, among them for the vector-valued
\cite{fewnomials} and complex analytic \cite{jlond} functions.

However, this method encounters serious difficulties when applied to
parametric families of functions, especially the families which admit
\emph{asymptotic} (rather than convergent) representation as in
\eqref{moura} as $t\to0^+$, with the exponents $\mu_j$ depending on
parameters and eventually coinciding between themselves. The operator
$\operD$, mentioned above, turns out to be very sensitive to these
parameters.

Our proof is based on introducing a new class of real
\emph{pseudo}differential operators $\operP$ which decrease in a
controlled way the number of isolated zeros. More precisely, for each
such operator $\operP$ and any family $\mathrm{F}$ of functions
analytic on the universal cover of the closed unit disk and
analytically depending on parameter there exists a finite number
(``Rolle index'') $N=N(\mathcal{P},\mathrm{F})$ such that for any
$u=u(t)\in\mathrm{F}$ on $(0,1)$ the numbers of real isolated zeros
of $u$ and $\mathcal{P}u$ are related by the ``Rolle inequality''
\begin{equation}\label{PR}
     \#\{u=0\}\le\#\{\operP u=0\}+N.
\end{equation}
These operators, depending on auxiliary parameters, in a sense
interpolate between the ``usual'' differential operators, the
Petrov difference operators that were used implicitly in
\cite{petrov} and explicitly in \cite{roitman-yakovenko}, and the
operator of ``taking imaginary part'' which played the key role in
\cite{quasialg}.

The operator $\operP$ is defined via the \emph{Mellin transform}
which associates with each function $u(t)$ (say, defined on the
interval $[0,1]$), the function $v(s)$ of a complex argument $s$ by
the formula
\begin{equation}\label{me}
    v=\M u,\qquad v(s)=\int_0^1t^{s-1}u(t)\,dt.
\end{equation}
We show (Theorem~\ref{lem:mellin} below) that the Mellin transform
$\M I$ of a pseudoabelian integral $I(t)$ is a function holomorphic
in the right half-plane $\Re s>s_0$, which admits analytic
continuation as a meromorphic function in the entire plane $s\in\C$
with poles located on finitely many real arithmetic progressions on
the shifted negative semiaxis $\R_-$. The distance between subsequent
poles may be very small, depending on the arithmetic of the exponents
$\l_j$, which results in the divergence of the representation
\eqref{moura}.

The operator $P$ is constructed as
\begin{equation}\label{p-init}
    P=\M^{-1}\operK\M,\qquad \operK:v(s)\mapsto K(s)v(s),
\end{equation}
where $K(s)$ is an entire function with zeros on these
progressions. Then the product $K\cdot \M I$ becomes an entire
function whose inverse Mellin transform $\M^{-1}$ is identically
zero. Note that a differential operator with constant coefficients
(the class sufficient for many applications of the
derivation--division algorithm) can be represented in this way
with a \emph{polynomial} kernel $K(s)$: this explains the
connections with the traditional method.

Flexibility of the construction of the operator $\operP$ (location of
zeros of $K(s)$) allows to choose $\operP=\operP_\theta$ in a way
\emph{analytically depending} on the form $\theta\in\mathscr D$ such
that $\operP_\theta I\equiv0$ for any pseudoabelian integral
$I=I_\theta(t)$ associated with  the integrable foliation $\theta=0$.

The cornerstone of the proof is the demonstration of the Rolle
inequality \eqref{PR} for the operator $\operP$. We prove that the
corresponding constant $N=N(\operP_\theta)$ depends in a uniformly
bounded way on the logarithmic 1-form $\theta\in\mathscr D$. This is
sufficient for the proof of Theorem~\ref{thm:main}.

\subsection{Acknowledgments} I grateful to S.
Yakovenko for drawing my attention to this problem and numerous
useful discussion. I am also grateful to C. Moura whose unpublished
preprint \cite{moura} was the starting point of all my
investigations. \pagebreak[4]
\section{Analytic representation of pseudoabelian
integrals}\label{sec:pai}

\subsection{Parametrization of admissible pairs of
forms}\label{ssec:parameters of pai} We define the complexification
$\setDC=\setDC(m_1,...,m_n)$ of $\setD$ as the space of logarithmic
forms $\theta$,
$$
\theta=\sum_{j=1}^n\l_j\frac{dp_j}{p_j},\qquad \l_j\in\C,\
p_j\in\C[x,y],\ \deg p_j\le m_j,
$$
where, as above, $p_j$ are  assumed to be coprime and irreducible.

Denote by  $\setD'$ the subset $\setD'\subset\setD\subset\setDC$ of
all $\theta\in\setD$ with smooth real algebraic curves $\{p_j=0\}$
intersecting transversally.

The pairs $(\theta,\omega)$ with $\theta\in \setD'$ and real
one-form $\omega$ admissible for $\theta$ with poles of order less
than $m_0$ form a set $\setB=\setB(m_0,m_1,...,m_n)$, which has an
evident topology of a real analytic vector bundle with base
$\setD'$.

The set $\setB$ is a real analytic subvariety of the set $\setBC$,
defined as the set of pairs of complex one-forms $(\theta, \omega)$,
with $\theta\in\setDC$ and $\omega$ admissible for $\theta$ with
poles of order less than $m_0$. The set $\setBC$ has an evident
structure of a complex analytic vector bundle over complex analytic
manifold $\setDC$.

\subsection{Local linearization of logarithmic foliations}
We start with a well-known fact whose demonstration is only briefly
sketched for completeness of exposition.

\begin{Lem}\label{lem:linearization works for c}
A closed meromorphic 1-form $\theta$ with first order poles on two
transversally intersecting analytic curves in $\C^2$ can be
analytically linearized near the point of their intersection. More
exactly, there is a linearizing transformation $\phi$ mapping a
neighborhood of a closed polydisk $\{|x|,|y|\le 2\}\subset\C^2$ to a
neighborhood of the point of the intersection, such that the
pull-back $\phi^*\theta$ is the logarithmic one-form $\l
\frac{dx}x+\mu \frac{dy}y$ defined in this neighborhood.

This linearizing transformation depends analytically on the complex
form $\theta\in\setDC$.
\end{Lem}

\begin{proof}
The two curves can be simultaneously rectified to become the
coordinate axes. The local fundamental group of the complement is
commutative, and is  generated by two loops around these axes: this
means that there exist two (complex) numbers $\l,\mu$ (the residues
of $\theta$) such that the difference $\theta-(\l \frac{dx}x+\mu
\frac{dy}y)$ is exact, the differential of a meromorphic function.
Because of the constraints on the order of the poles, this difference
is the differential of a holomorphic function $dg$, holomorphically
depending on $\theta$. The inverse of the holomorphic transformation
$(x,y)\mapsto (x,y\exp(\mu^{-1}g))$ brings the form $\theta$ to the
required ``linear'' normal form in some neighborhood of the origin,
which can be subsequently expanded to the polydisk  by a linear
dilatation in $x,y$.
\end{proof}

Here is the linearization transformation written explicitly in terms
of the first integral. Let
\begin{equation}\label{eq:FI}
f=Cp_1(x,y)^{\l}p_2(x,y)^{\mu}\prod_{j=3}^n p_j(x,y)^{\lambda_j}
\end{equation}
be a first integral of $\theta$, $C\neq0$, and assume that the
transversally intersecting curves are $\{p_1=0\}$ and $\{p_2=0\}$.
The inverse of the mapping
$$(x,y)\to (\tilde{x},\tilde{y}),\quad\text{where}\
\tilde{x}=p_1(x,y), \  \tilde{y}=p_2\prod_{j=3}^n
p_j(x,y)^{\lambda_j/\mu}$$ is the required diffeomorphism $\phi$.

\begin{Lem}\label{lem:first integral is a monomial}
For any  $C$ bigger than some $C_0\gg1$ the dilatation can be chosen
in such a way that the first integral \eqref{eq:FI} becomes a
monomial $\tilde{x}^{\lambda}\tilde{y}^{\mu}$ in the new
coordinates.\end{Lem}

Indeed, after linearization as in Lemma \ref{lem:linearization works
for c} the first integral becomes
$\const\tilde{x}^{\lambda}\tilde{y}^{\mu}$, and if $\const>1$, then
an additional dilatation brings the first integral to the required
form.\qed

\subsection{Representation for pseudoabelian integral}

\subsubsection{Local computation near saddles}
We compute $I(t)$ by evaluating separately integrals along pieces
$\delta_{ij}$ of $L_t$ lying near saddles $\{p_i=p_j=0\}$ and along
the pieces $\delta_j$ of $L_t$ lying near smooth pieces of
$\{p_j=0\}$.

\begin{figure}
\centerline{\epsfysize=0.25\vsize\epsffile{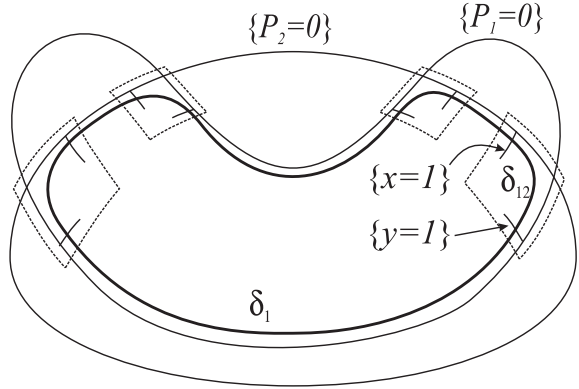}}
\caption{Splitting of a cycle into pieces.}\label{fig:cycle}
\end{figure}

Integrals of monomial 1-forms along the piece $\delta_{ij}$ of the
integral curve $L_t=\{f=t\}$ lying between two cross-sections,
$\{x=1\}$  and $\{y=1\}$, can be explicitly computed. Indeed, in the
linearized coordinates $(x,y)$ the curve $\delta_{ij}$  is, by
Lemma~\ref{lem:first integral is a monomial}, given by the explicit
formula
\begin{equation}\label{corner}
   \delta_{ij}=\{y=t^{1/\mu}x^{-\l/\mu}\},\qquad t\in(0,1),
    x\in(t^{1/\l},1),
\end{equation}
eventually after replacing $f$ by $Cf$, with $C\gg 1$. Therefore the
integration of $\omega=x^{p-1}y^q\,dx$, $(p,q)\in\mathbb Z^2$, is
simple:
\begin{equation}\label{I-mono}
    \int_{t^{1/\l}}^1x^{p-1}y^q\,dx=
    t^{q/\mu}\,\frac{x^{p-q\l/\mu}}{p-q\l/\mu}\bigg|_{t^{1/\l}}^1=
  \l^{-1}\frac{t^{q/\mu}-t^{p/\l}}{p/\l-q/\mu},
\end{equation}
if $q/\mu\ne p/\l$. For the resonant case $q/\mu=p/\l$ the same
computation yields the answer $-\l^{-1}t^{q/\mu}\log t$ which can be
otherwise obtained by direct passage to limit.  Computation of
$\int_{\delta_{ij}}x^py^{q-1}\,dy$ is similar.

Denote by $\ell_{pq\l\mu}(t)$ the following elementary function:
\begin{equation}
 \ell_{pq\l\mu}(t)=
    \begin{cases}
    \dfrac{t^{p/\l}-t^{q/\mu}}{p/\l-q/\mu},&\text{if }p/\l\ne q/\mu,
    \\[10pt]
    {t^{p/\l}\log t}&\text{otherwise},
    \end{cases}
\end{equation}
One can verify that $\ell_{pq\l\mu}(t)$ depends analytically on
$\lambda,\mu$ for any $t>0$.

\subsubsection{Analytic representation of pseudoabelian integral}

\begin{Lem}\label{lem:pai}
After some rescaling of $t$, the pseudoabelian integral $I(t)$
corresponding to $(\theta,\omega)$ sufficiently close to
$(\theta_0,\omega_0)$ in $\mathscr{B}$   admits the representation
near the point $t=0$
\begin{equation}\label{eq:pai}
    I(t)=\sum_{p,q,\l,\mu}a_{pq\l\mu}\ell_{pq\l\mu}
    (t)+\sum_{r,\l}b_{r\l}t^{r/\l},
\end{equation}
where  $p,q,r$ run over a translated octant $\mathbb Z_+^3$,
$p,q,r\in -m_0+\N$, and the real indices $\l,\mu$ run independently
over the set $\L=\{\l_1,\dots,\l_n\}\subset\R^n_+$ of the Darboux
exponents of the closed 1-form $\theta$.

Moreover, the coefficients $a_{pq\l\mu},b_{r\l}$ depend analytically
on the forms $\theta,\omega$, and can be continued holomorphically to
some neighborhood $U$ of $(\theta_0,\omega_0)\in\setBC$, and
continuously to its closure $\overline{U}\Subset\setBC$. The
coefficients $a_{\bullet\l\mu},b_{\bullet\l}$ decrease at least
exponentially: there exist $C,\rho>2$ such that the following bounds
hold in $\overline{U}$
\begin{equation}\label{eq:pai bounds}
    |a_{pq\l\mu}|\le C\rho^{-(p+q)}, |b_{r\l}|\le
    C\rho^{-r}.
\end{equation}
\end{Lem}

\begin{proof}
By assumptions, the Hausdorff limit $\overline
L_0=\lim_{t\to0^-}L_t\subset S(\theta)$ of the real ovals is a
bounded curvilinear separatrix polygon, whose edges and vertices
depend analytically on $\theta\in\setD$ in a sufficiently small
neighborhood of $\theta_0$ by the Implicit Function Theorem.

Near each vertex one can linearize the foliation using Corollary
\ref{lem:linearization works for c} in such a way that the first
integral is a monomial in the new coordinates,
$f={x}^{\lambda}{y}^{\mu}$, eventually after a rescaling of $f$.

The linearizing transformation brings the  rational form $\omega$
into a meromorphic form in the bidisk $\{|x|,|y|\le 2\}$ with poles
of order less than $m_0$ on the coordinate axes. Such form can be
represented as a series
$$
\sum_{p,q> -m_0} \alpha_{pq}x^{p-1}y^q\,dx+\beta_{pq}x^py^{q-1}\,dy
$$
converging in the bidisk $\{|x|,|y|\le 2\}$. This convergence implies
the exponentially decreasing upper bounds \eqref{eq:pai bounds} on
the Taylor coefficients $\alpha_{pq},\beta_{pq}$ of the form $\omega$
in the linearizing coordinates $x,y$. Since the linearization depends
analytically on $\theta$ in a closure $\bar{U}\Subset\setBC$ of a
sufficiently small neighborhood $U$ of $\theta_0\in\setBC$ by
Lemma~\ref{lem:linearization works for c}, the coefficients
$\alpha_{pq},\beta{pq}$ are analytic in $\overline{U}$ and the upper
bounds are uniform in this neighborhood. Integrating the series
termwise along $\delta_{ij}(t)=\{f=t\}\cap\{|x|,|y|\le 2\}$ as in
\eqref{I-mono}, and summing up over all vertices of $\overline{L}_0$
corresponding to the pair $(\l,\mu)$ (their number is finite by
Bezout theorem) we obtain the first sum in \eqref{eq:pai} and the
exponential upper bounds \eqref{eq:pai bounds} for $a_{pq\l\mu}$, as
finite linear combinations of $\alpha_{pq}, \beta_{pq}$.

The local cross-sections near the corners given in linearized
coordinates by $\{x=1\}$ and $\{y=1\}$ become analytic curves
transversal to $S(\theta)$ in the initial coordinates on the plane.
The integral of the form $\omega$ along arcs of leaves $L_t$ between
the two cross-sections near two endpoints of the $j$th edge (on the
curve $\{p_j=0\}$), depends analytically on the arc, as the latter
remains nonsingular. In the chart $t$ obtained by restriction of the
Darbouxian integral, this means that the corresponding contribution
to the integral is an analytic function of $t^{1/\l_j}$, where $\l_j$
is the residue (Darboux exponent) corresponding to the edge. Again,
after a rescaling one can assume that this analytic function
converges in a small neighborhood of $\{|t^{1/\l_j}|\le 2\}$,
uniformly over a small neighborhood of
 $(\theta_0,\omega_0)\in\setBC$, thus giving the second inequality in
\eqref{eq:pai bounds}.
\end{proof}

\section{The class $\spaceJ$ and its Mellin transform}\label{sec:J}
We introduce a class $\spaceJ$  of functions of one variable as sums
of series of type \eqref{eq:pai} with coefficients satisfying
\eqref{eq:pai bounds}, and define the notion of analytic
$\spaceJ$-family in Section~\ref{sec:J fam and E} using the statement
of Lemma~\ref{lem:pai} as its definition. We prove in
Section~\ref{sec:petrov} that the number of zeros on $[0,1]$ is
uniformly bounded for any analytic $\spaceJ$-family.

\subsection{ Definition of the class $\spaceJ$}\label{ssec: def J}
\begin{Def}\label{def:J}

We define  $\spaceJ$-series as a formal series $\sigma$ of the form
\begin{equation}\label{eq:def f in J}
\sum_{p,q,\l,\mu}a_{pq\l\mu}\ell_{pq\l\mu}(t)+\sum_{r,\l}b_{r\l}t^{r/\l},
\qquad a_{pq\l\mu},b_{r\l}\in\C,
\end{equation}
where $\l,\mu\in\Lambda=\{\l_1,...,\l_n\}\subset\R_+$ - a finite set,
 $p,q,r\in -m_\sigma+\N$,  and
\begin{equation}\label{eq:def f in J bounds}
|a_{pq\l\mu}|<C\rho^{-p-q}, |b_{r\l}|<C\rho^{-r}\quad\text{for some
}C,\rho>2.\end{equation}

The set $\Lambda$ is called the spectrum of $\sigma$.

We denote further by $M=M(\sigma)$ some number smaller than all
exponents in \eqref{eq:def f in J}, e.g.:
\begin{equation}\label{eq: def of M}
M=\min_i(-m_\sigma\lambda_i).
\end{equation}

The class $\spaceJ$ is the class of functions $f(t)$ in one variable
$t\in[0,1]$, called $\spaceJ$-functions, representable as a sum of
some $\spaceJ$-series.
\end{Def}

One should emphasize that the $\spaceJ$-series corresponding to a
function $f\in\spaceJ$ is not unique: indeed, even the function
$\ell_{pq\l\mu}(t)$ can be written as a sum of two monomials. Also,
we do not know any effective way to check if a given function belongs
to $\spaceJ$. However, according to Lemma \ref{lem:pai}, the
pseudoabelian integrals belong to $\spaceJ$.

\subsection{Analytic continuation of
$\spaceJ$-functions}\label{ssec:cont of J} We prove here that
$\spaceJ$-functions can be analytically continued to the universal
cover of the punctured unit disk.

We start with an elementary inequality.
\begin{Lem}\label{lem:bound for compensator}
For any $t,\alpha,\beta\in\C$ we have
$$
\left|\frac{t^{\alpha}-t^{\beta}}{\alpha-\beta}\right|\le|t^\gamma||\log
t|
$$
for some $\gamma\in[\alpha,\beta]\subset\C$.
\end{Lem}

Indeed, take $z=\log t$ in the following inequality
$$
\left|e^{\alpha z}-e^{\beta
z}\right|=\left|z\int_{\beta}^{\alpha}e^{wz}dw\right|\le|\alpha-\beta|\cdot\left|e^{\gamma
z}\right|\cdot|z|\quad \text{for some } \gamma\in[\beta,\alpha].\qed
$$

Let $\widetilde{\C^*}$ denote the universal cover of
$\C^*=\C\setminus\{0\}$ and let $\tilde
D^\circ\subset\widetilde{\C^*}$ be the universal cover of the
punctured closed unit disk.
\begin{Lem}\label{lem:uniform convergence}
For any $A>0$ the series \eqref{eq:def f in J} converges uniformly
after multiplication by $t^{-M}$ in a neighborhood in
$\widetilde{\C^*}$ of the sector $S_A=\{|t|\le 1, |\arg t|\le
A\}\subset\tilde D^\circ$, and, in particular, on $[0,1]$.
\end{Lem}

Lemma \ref{lem:bound for compensator} implies that for some constant
$C_A$ and for  all $p,q,\l,\mu$ in \eqref{eq:def f in J} we have
$|t^{-M} \ell_{pq\l\mu}(t)|\le C_A 2^{(p+q)}$ as long as
$|t|^{1/\l},|t|^{1/\mu}<2$.  These inequalities define a neighborhood
of $S_A$ in which the series \eqref{eq:def f in J} converges
uniformly by \eqref{eq:def f in J bounds}.\qed

\begin{Cor}
$\spaceJ$-functions are real analytic functions on $(0,1]$ admitting
analytic continuation to the universal cover of the punctured closed
unit disc.
\end{Cor}

\begin{small}
\begin{Rem}
If the ratios of all consecutive Darboux exponents $\l_i/\l_{j}$,
$i,j=1,\dots,n$ are either rational, or ``nice'' irrational (not
admitting too close approximations by rational numbers so that the
differences $p\l_i-q\l_j$ decrease no faster than exponentially with
$p,q\to+\infty$), then the series \eqref{eq:pai} can be re-expanded
in powers of $t^{1/\l_j}$, and this re-expansion \emph{converges}.
Then the integral $I(t)$ becomes representable under the form
\eqref{moura}. Application of the standard fewnomial technique
allowed C.~Moura to prove  that for such collections of the Darboux
exponents the pseudoabelian integral $I(t)$ is non-oscillating,
moreover, that the number of isolated zeros is locally uniformly
bounded over $\omega$, \emph{provided that the Darbouxian exponents
are fixed} \cite{moura}. However, variation of the exponents clearly
destroys the arithmetic conditions.
\end{Rem}

\begin{Rem}
The transcendental binomials $\ell_{pq\l\mu}$ which in the limit
converge to logarithms, first appeared at least as early as in 1951
in the work by E.~Leontovich \cite{leontovich} on bifurcations of
separatrix loop. Since then they regularly re-appear in different
disguises, most recently as the so called
Ecalle--Roussarie compensator.

This ubiquity can be probably explained by the fact that the
function $\ell_\e(t)=t(t^\e-1)/\e$ is in a sense the only possible
Pfaffian deformation of the ``logarithm'' $\ell_0(t)=-t\log t$.
Indeed, the graph of the latter function is an integral curve of the
linear equation $\xi_0=0$, where $\xi_0=t\,dy+(y-t)\,dt$ is a
Pfaffian form. This linear Pfaffian equation has a singularity of
the Poincar\'e type (resonant node) at the origin. Any (nonlinear)
deformation $\xi_\e$ of the form $\xi_0$ is necessarily analytically
linearizable by the Poincar\'e theorem, the linearizing chart
depending analytically on the parameter $\e$. Among integral curves
of the perturbed linear system one is the graph of $\ell_\e$.
\end{Rem}
\end{small}

\subsection{Mellin transform of $\spaceJ$}\label{ssec:mellin}

The (one-sided) Mellin transform of a function $u(t)\in
L^1_{loc}((0,1])$ defined on the interval $(0,1]$ is the function
$v=\M u$ of complex variable $s$ defined by the integral transform
\eqref{me}:
$$
v(s)=\M u(s)=\int_0^1t^{s-1}u(t)dt.
$$

If the function $u$ grows moderately as $t\to 0^+$,
$|u(t)|=o(t^{-M})$ for some finite $M$, then the integral in
\eqref{me} converges uniformly on compact subsets of the half-plane
$\{\Re s> M\}$. In particular, Mellin transforms are defined for all
$\spaceJ$-functions.

The formal Mellin transform of a $\spaceJ$-series is defined by its
action on the basic functions
\begin{equation*}
\M(t^{r/\l})=\frac 1 {s+r/\l},\quad \M (\ell_{pq\l\mu})=\frac 1
{(s+p/\l)(s+q/\mu)}
\end{equation*}
and extended by linearity, so the result is a formal sum of these
rational functions. We prove that this formal sum converges to a
meromorphic function which is an analytic continuation of the Mellin
transform of the sum of the $\spaceJ$-series.

\begin{Lem}\label{lem:mellin}
The Mellin transform $g=\M f$ of a function $f\in \spaceJ$ given as a
sum of the series \eqref{eq:def f in J} is the sum of the series:
\begin{equation}\label{eq:g}
g(s)=\sum_{p,q,\l,\mu}\frac{a_{pq\l\mu}}{(s+p/\l)(s+q/\mu)}+
\sum_{r,\l}\frac{b_{r\l}}{s+r/\l},
\end{equation}
with the same $p,q,r,\l,\mu$ and $a_{pq\l\mu},b_{r\l}$ as in
Definition \ref{def:J}.

In particular,  $g(s)$ extends analytically as a meromorphic function
on the entire complex plane $\C$, the poles of $g(s)$ are of at most
second order, are all real and lie in the union of $n$ arithmetic
progressions of the form $\cup_{i=1}^n \lambda_i^{-1}(m_0-\N)$.
\end{Lem}

\begin{proof}
The Mellin transforms of each elementary binomial $\ell_{pq\l\mu}$
and of $t^{r/\l}$ are equal to its formal Mellin transform. Termwise
integration of right-hand side of \eqref{eq:def f in J} multiplied by
$t^{s-1}$ is possible due to the uniform convergence on $[0,1]$ of
this product  for $\Re s>-M$ by Corollary \ref{lem:uniform
convergence}. Therefore we get \eqref{eq:g} in $\{\Re s>-M\}$. This
double series converges to a meromorphic function on $\C$ given by
the same formula.
\end{proof}

Let $\gamma_M$ be the boundary of the semistrip $\{z|\Re z\le
-M+1,|\Im z|\le -1\}$ oriented counterclockwise, with $M$ defined in
\eqref{eq: def of M}.

For any function bounded on $\gamma_M$ we define an integral
transform $\operC_M$ as:
\begin{equation}\label{eq:invMellin}
\operC_Mg(t)=\frac 1 {2\pi i}\int_\gamma t^{-s}g(s)ds.
\end{equation}

For a meromorphic in $\C$ function $g$ bounded on $\{|\Im s|=1\}$
with poles on the real line and bounded from the right we define
$\operC g$ as $\lim_{M\to-\infty}\operC_M g$.

\begin{Thm}\label{thm:inv mellin}
The Mellin transform $\M f$ of any function $f\in\spaceJ$ is bounded
on $\gamma_M$, where $M=M(\sigma)$ for some $\spaceJ$-series
representing $f$.

The restriction of the integral transform  $\operC$ to $\M\spaceJ$ is
the inverse operator to the Mellin transform $\M$.

For $\kappa\in \R$ and $f\in\spaceJ$ the transform $\operC$ of
$e^{-i\kappa s}\M f(s)$  is defined and is equal to
$f(e^{i\kappa}t)$\/\emph{:}
\begin{equation}\label{eq:shift}
\operC(e^{-i\kappa s}\M f(s))=f(e^{i\kappa}t).
\end{equation}
\end{Thm}

Abusing the language we call $\operC$ the inverse Mellin transform,
though its domain of definition is bigger than $\M\spaceJ$.

\begin{proof}
Evidently $|s+p/\l|\ge 1$ on $\gamma$, so the convergence of
\eqref{eq:g} is uniform and $|g|\le
\sum_{p,q,\l,\mu}|{a_{pq\l\mu}}|+\sum_{r,\l}|{b_{r\l}}|<\infty$ on
$\gamma$.

The kernel $t^{-s}$ decreases exponentially on $\gamma$ for $|t|<1$,
so we can integrate $t^{-s}g(s)$ termwise on $\gamma$. Calculating
the residues, we see that
\begin{equation*}
\operC\left(\frac 1
{s+r/\l}\right)=t^{r/\l},\quad\operC\left(\frac 1
{(s+p/\l)(s+q/\mu)}\right)=\ell_{pq\l\mu},
\end{equation*}
so the Mellin transform and the inverse Mellin transform translate
\eqref{eq:def f in J} to \eqref{eq:g} and vice versa, i.e. are
mutually inverse if restricted to $\spaceJ$ and $\M\spaceJ$
correspondingly.

Since $e^{-i\kappa s}$ is bounded by $e^{|\kappa|}$ on $\gamma$, one
can compute the inverse Mellin transform \eqref{eq:invMellin}
termwise, again due to the uniform convergence of the series for
$g(s)$ and the exponential decay of $t^{-s}$. For $\ell_{pq\l\mu}$
and $t^\alpha$ the \eqref{eq:shift} can be checked immediately.

\end{proof}

\subsection{Asymptotic series for $\spaceJ$-functions and
quasianalyticity of $\spaceJ$}\label{ssec:quasianalyticity}

Let $f\in \spaceJ$ be a sum of some $\spaceJ$-series $\sigma$. Any
$\spaceJ$-series $\sigma$ can be expanded into a formal series by
powers of $t$:
\begin{equation}\label{eq:asympt ser for f}
\check{f}=\sum_{\alpha} t^{\alpha}\left(c^{-1}_\alpha+c^{-2}_\alpha\log
t\right),\qquad \alpha\in (-m_\sigma+\N)\Lambda^{-1},
\end{equation}
where  $\Lambda$ is the spectrum of $\sigma$ and $\Lambda^{-1}$
denotes the set of its reciprocals. We prove below that the result of
the expansion does not depend on the choice of $\sigma$ and
represents the asymptotic series for $f$. Moreover, the class
$\spaceJ$ turns out to be quasianalytic at $0$: any two functions
having the same asymptotic series coincide.

\begin{Lem}\label{lem:asympt for f}
The coefficients $c^{-1}_\alpha,c^{-2}_\alpha$ are the Laurent
coefficients of the Mellin transform $g=\M f$ of $f$ at the pole
$s=\alpha$.

The series $\check f$ is an asymptotic series for $f$, so does not
depend on $\sigma$. If $\check f\equiv0$ then $f\equiv 0$.
Equivalently, the only entire function in $\M\spaceJ$ is zero.
\end{Lem}
\begin{proof}
The first claim follows because the Laurent coefficients of $g$ are
given by exactly the same formulae as $c_\alpha,d_\alpha$:
$$
c^{-1}_\alpha=\sum_{\substack{q,\mu\\q/\mu=\alpha}}b_{q\mu}+
\sum_{\substack{p,q,\l,\mu\\q/\mu\neq p/\l=\alpha}}\frac{a_{pq\l\mu}}{p/\l-q/\mu},
\quad c^{-2}_\alpha=\sum_{\substack{p,q,\l\mu\\q/\mu=p/\l=\alpha}}a_{pq\l\mu},
$$
both converging due to \eqref{eq:def f in J bounds}.

Choose some $0<\l_0\le \min \Lambda$. Each segment $[k\l_0,
(k+1)\l_0)$, $k\in\Z$, has at most $n$ points common with
$\Lambda\Z$, so one can choose $M_k$ from this segment at least
$\l_0/2n$-distant from $\Lambda\Z$.

\begin{figure}\label{fig:asympt}
\centerline{\epsfysize=0.20\vsize\epsffile{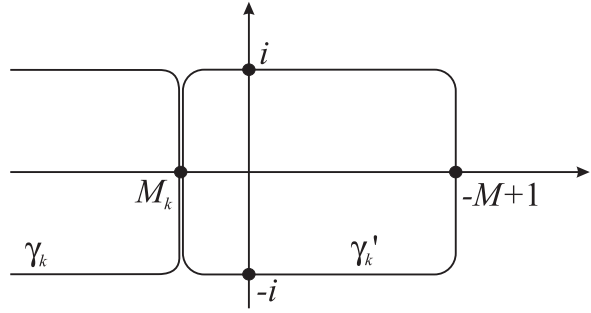}}
\caption{Splitting of the integration contour $\gamma$ into two.}
\end{figure}

This implies that $|z+p/\l|\ge\l_0/2n$ on $\Re z=M_k$, so $|g|\le C$
on the boundary $\gamma_k$ of the semistrip  $\{\Re z\le M_k, |\Im
z|\le 1\}$, and the main point is that $C$ is independent of $k$.
Denote by $\gamma'_k$ the boundary of the rectangular $\{M_k\le\Re
z\le -M+1, |\Im z|\le 1\}$ and let split integration along $\gamma$
in the definition of $\operC$ into a sum of two integrals, one along
$\gamma_k$ and another along $\gamma'_k$, see
Figure~\ref{fig:asympt}.

Evaluation of the integral along $\gamma'_k$ gives a finite partial
sum of $\check f$:
$$
\frac 1 {2\pi i}\int_{\gamma'_k} t^{-s}g(s)ds=\sum_{\alpha< -M_k}
t^{\alpha}\left(c_\alpha+d_\alpha\log t\right).
$$

Integral along $\gamma_k$ can be estimated from above by
$2Ct^{-M_k}\left(1+1/\log t\right)$, which proves the second claim.

If $\check f\equiv0$, then $g$ has no poles at all, i.e. $g$ is
entire. Therefore the integral along $\gamma'_k$ is zero for all $k$.
Therefore $|f(t)|\le 2Ct^{-M_k}\left(1+1/\log t\right)$ for all
$t\in(0,1)$, and for all $M_k$. Since $M_k\to-\infty$ as
$k\to-\infty$ and $C$ is independent of $k$, we conclude that
$f(t)\equiv0$ on $(0,1)$, and therefore everywhere by analyticity.
This means that $g=\M f\equiv 0$ as well.\end{proof}

\section{Analytic $\M \spaceJ$-families and operator
$\operE_{\kappa}$}\label{sec:J fam and E}
\begin{Def}
Let $V\Subset\R^n$ be some compact set, and let $U\subset\C^n$ be a
bounded neighborhood of $V$. We denote by $\spaceO_{V\subset U}$ the
space of real analytic functions on $V$ which can be extended
holomorphically to $U$ and continuously to the closure $\overline{U}$
of $U$, equipped with the norm $\| u\|_{\spaceO}:=
\max_{\nu\in\overline{U}}|u(\nu|)$.
\end{Def}
\begin{Lem}\label{lem:O complete}
The space $\spaceO_{V\subset U}$ is complete.\qed \end{Lem}

\begin{Def}\label{def:an fam in J}
Consider a family $\check{\mathrm{F}}=\{\sigma_\nu\}$ of
$\spaceJ$-series, $\nu\in V\Subset\R^N$, and denote by
$\Lambda(\nu)=\{\l_1(\nu),...,\l_n(\nu)\}\subset\R_+$ the spectrum of
the $\spaceJ$-series $\sigma(\nu)$.

We say that $\check{\mathrm{F}}$ is an analytic family of
$\spaceJ$-series on pair $V\subset U$ for some bounded neighborhood
$U\subset\C^N$ of $V$ in $\C^N$ if the functions $\l_i(\nu)$,
$a_{pq\l(\nu)\mu(\nu)}(\nu)$, $b_{r\l(\nu)}(\nu)$ belong to
$\spaceO_{V\subset U}$, the analytic continuation of $\l_i(\nu)$ does
not vanish in $\overline{U}$ (i.e.
$\l_i^{-1}(\nu)\in\spaceO_{V\subset U}$ as well), and
\begin{equation}\label{eq:coeff bounded in family}
\left\|a_{pq\l(\nu)\mu(\nu)}(\nu)\right\|_{\spaceO}|\le C\rho^{-(p+q)},\quad
\left\|b_{r\l(\nu)}(\nu)\right\|_{\spaceO}\le C\rho^{-r}
\end{equation}
for some $C=C(\check{\mathrm{F}}),\rho=\rho(\check{\mathrm{F}})>2$.

An analytic $\spaceJ$-family $\mathrm{F}=\{f(t;\nu)\}\subset\spaceJ$
is a family of sums of an  analytic family of $\spaceJ$-series. In
other words, for every $f(t;\nu)$ one  can choose a $\spaceJ$-series
$\sigma_\nu$ such that $f(t;\nu)$ is a sum of $\sigma_\nu$, and
$\{\sigma_\nu\}$ is an  analytic family of $\spaceJ$-series.

Analytic $\M \spaceJ$-family is by definition the Mellin transform of
an analytic $\spaceJ$-family.
\end{Def}

\begin{Lem}\label{lem:m is unif bdd in fam}
For a given  analytic family $\check{\mathrm{F}}$ of $\spaceJ$-series
the number $m_{\sigma}$ is uniformly bounded from above by some
number $m_{\mathrm{F}}<\infty$ for all
$\sigma=\sigma_\nu\in\check{\mathrm{F}}$.
\end{Lem}
Indeed, the sets $S_k=\{\nu\in U| m_{\sigma_\nu}\le k\in\N\}$ are
analytic subsets of $U$ and their union is $U$, so $S_k=U$ for some
$k$.\qed

Lemma \ref{lem:pai} says that the family of pseudoabelian integrals
parameterized by $\mathscr{B}$ is locally an analytic
$\spaceJ$-family.

Let $\check{\mathrm{F}}=\{\sigma_\nu\}$ be an analytic family of
$\spaceJ$-series on pair $V\subset U$, and consider the series
\begin{equation}\label{eq:continuation family}
\sum_{p,q,\l(\nu),\mu(\nu)}a_{pq\l(\nu)\mu(\nu)}\ell_{pq\l(\nu)\mu(\nu)}(t)+
\sum_{r,\l(\nu)}b_{r\l(\nu)}t^{r/\l(\nu)}\quad\text{for}\nu\in U.
\end{equation}

\begin{Lem}\label{lem:continuation family}
For every $A>0$ there exists some $\varepsilon>0$ such that this
series converges uniformly on compact subsets of
$$U_\varepsilon=\{|\arg t|\le
A, 0<|t|< e^\varepsilon\}\times\(\{|\Im \l^{-1}_i|< \varepsilon\}\cap
U\)\subset\widetilde{\C^*}\times U,$$ so its sum is a function
holomorphic on $U_\varepsilon$.
\end{Lem}

\begin{proof}
In $U_\varepsilon$  we get $|t^{\l^{-1}(\nu)}|\le
\exp\left(\varepsilon\left[A+\|\l^{-1}(\nu)\|_{\spaceO}\right]\right)$,
i.e. $|\ell_{pq\l\mu}(t)|$ can grow exponentially fast. However, the
fast decay of coefficients $a_{pq\l\mu},b_{r\l}$ guarantees
convergence if $\varepsilon$ is sufficiently small: the series
\eqref{eq:def f in J} converges uniformly on compact subsets of
$U_\varepsilon$ as soon as
$\varepsilon\left[A+\|\l_i^{-1}(\nu)\|_{\spaceO}\right]<\ln 2$ for
$i=1,...,n$.\end{proof}

\begin{Rem} For functions from $\M \spaceJ$ one can try to define analytic
dependence on parameter as analytic dependence of $f(t;\nu)$ on $\nu$
on compact subsets of $\C\smallsetminus\R$. Analytic $\M
\spaceJ$-families depend analytically on parameter in this sense, as
shown above, but the opposite implication seems to be wrong.\end{Rem}

\subsection{The operator $\operE_{\kappa}$}\label{ssec:operE}
From Lemma~\ref{lem:e is well-defined} below follows that the
rotation of argument preserves $\spaceJ$: for any $f(t)\in\spaceJ$
the function $f(e^{i\kappa}t)$, $\kappa\in\R$, is also in $\spaceJ$.
From \eqref{eq:shift} it seems that the Mellin conjugate to the
rotation is the operator of multiplication by $e^{-i\kappa s}$.
However, the latter does not preserves $\M\spaceJ$: the function
$e^{-i\kappa s}\M f(s)$ grows as $|\Im s|\to\infty$, while any
function in $\M\spaceJ$ decreases as $|\Im s|^{-1}$.

In this section we define the Mellin conjugate $\operE_\kappa$ of the
rotation on $\M\spaceJ$, and prove that application of operator
$\operE_\kappa$ with $\kappa$ analytically depending on parameters
preserves $\M \spaceJ$-families.

Let start from the standard computation applied to $f\in\spaceJ$:
\begin{multline}\label{eq:rotation}
\M(f(e^{i\kappa}t))=\int_0^1t^{s-1}f(e^{i\kappa}t)dt\\
=\int_0^{e^{-i\kappa}}t^{s-1}f(e^{i\kappa}t)dt+
\int_{C_\kappa}t^{s-1}f(e^{i\kappa}t)dt \qquad\qquad\qquad\\
=e^{-i\kappa s}\int_0^{1}u^{s-1}f(u)du+R(s)=e^{-i\kappa s}\M
f(s)+R(s),
\end{multline}
where $C_\kappa=\{e^{iz},z\in[-\kappa,0]\}$ and
$R(s)=\int_{C_{\kappa}}t^{s-1}f(t)dt$ is an entire function of
exponential type. One can immediately see that $R(s)$ is bounded in
$\{|\Im s|\le 1\}$. Therefore $\operC R\equiv 0$, which explains
\eqref{eq:shift}. We take this as a definition of $\operE_\kappa$. In
a sense $\operE_\kappa$ is an operator of projection of $e^{-i\kappa
s}\M\spaceJ$ to $\M\spaceJ$ along the space of entire functions.

\begin{Def}\label{def:ekappa}
Let $\kappa\in \R$. For any $g\in\M\spaceJ$  define $\operE_\kappa g$
as the only function in $\M \spaceJ$ such that the difference
$e^{-i\kappa s}g-\operE_\kappa(g)$ is an entire function.\end{Def}

\begin{Rem}
One can show that  $\operE_{\kappa}g(s)$ is an analytic continuation
from $\{\Im s>1\}$ of the integral
$$
\frac 1 {2\pi i}\int_\gamma\frac{e^{-i\kappa z}g}{z-s}dz,
$$
where $\gamma$ is as in \eqref{eq:invMellin}.
\end{Rem}

\begin{Lem}\label{lem:e is well-defined}
Let $\mathrm{F}=\{g(s;\nu)\}$ be an analytic $\M \spaceJ$-family on
pair $V\subset U$, and let $\kappa=\kappa(\nu)\in\spaceO_{V\subset
U}$. Then
\begin{enumerate}
\item $\{\operE_{\kappa(\nu)}g(s;\nu)\}$ is defined and is an
    analytic $\M \spaceJ$-family on pair $V\subset \widetilde{U}$
    for some $\widetilde{U}\subset U$.
\item the operator $\operC\circ \operE_{\kappa}\circ \M$ maps
    analytic $\spaceJ$-families to analytic $\spaceJ$-families,
\item $\operC\operE_\kappa \M f=(e^{i\kappa}t)$, and
\item the operator of rotation of argument preserves $\spaceJ$.
\end{enumerate}
\end{Lem}

\begin{proof} The uniqueness in Definition~\ref{def:ekappa} follows
from Lemma~\ref{lem:asympt for f}: let $g_1,g_2\in\M \spaceJ$ be  two
functions which differ from $e^{-i\kappa s}g$ by entire functions.
Then their difference is entire and lies in $\M \spaceJ$, which, by
Lemma \ref{lem:asympt for f}, implies that the difference is zero.

Let us construct $\operE_\kappa$. This is simple for $\frac 1
{(s+p/\l)(s+q/\mu)}$ and $\frac 1 {s+r/\l}$:
$$
\frac{e^{-i\kappa s}}{s-x}=\frac{e^{-i\kappa
x}}{s-x}+\frac{e^{-i\kappa s}-e^{-i\kappa x}}{s-x},
$$
and
\begin{equation*}
2\frac{e^{-i\kappa s}}{(s-x)(s-y)}=\frac{e^{-i\kappa x}+e^{-i\kappa
y}}{(s-x)(s-y)}  +\frac{e^{-i\kappa x}-e^{-i\kappa y}}{x-y}\(\frac
1{s-x}+\frac1{s-y}\)+R,
\end{equation*}
where $R=2\frac{(s-x)e^{-i\kappa y}+(x-y)e^{-i\kappa
s}+(y-s)e^{-i\kappa x}}{(s-x)(s-y)(x-y)}$ is entire in $s$.

We extend this by linearity to any $g\in\M\spaceJ$. Namely, for $g$
given by \eqref{eq:g} we define $\operE_\kappa g(s)$
$$
\operE_\kappa
g(s)=\sum_{p,q,\l,\mu}\frac{\tilde{a}_{pq\l\mu}}{(s+p/\l)(s+q/\mu)}+
\sum_{r,\l}\frac{\tilde{b}_{r\l}}{s+r/\l},
$$
by the following formulae:
\begin{gather}
\tilde{a}_{pq\l\mu}=\frac 1
2\(e^{-i\kappa p/\lambda}+e^{-i\kappa q/\mu}\)a_{pq\l\mu}, \label{eq:acute a}\\
\tilde{b}_{r\l}=e^{-i\kappa r/\lambda}b_{r\l}
+\sum_{q,\mu}\frac{e^{-i\kappa r/\lambda}-e^{-i\kappa q/\mu}}
    {r/\lambda-q/\mu}a_{rq\l\mu}.\label{eq:acute b}
\end{gather}

Take $\varepsilon$ so small that $\tilde{\rho}=\rho
e^{-\varepsilon}>2$. Let a neighborhood $\widetilde{U}\subset U$
of $V$ be so small that $\left|\Im
\left(\kappa(\nu)\lambda_i^{-1}(\nu)\right)\right|<\varepsilon$
for $i=1,...,n$ and  $\nu\in\tilde{U}$. Then $|e^{-i\kappa
r/\lambda}|<e^{\varepsilon r}$ in $\widetilde{U}$.

For $p,q,r>-m$ as in \eqref{eq:g} we have $\max (p,q)<p+q+m$, so
for $\nu\in\widetilde{U}$ we have estimates
$$
\frac 1 2\left|e^{-i\kappa p/\lambda}+e^{-i\kappa q/\mu}\right|\le
e^{\varepsilon(p+q+m)}, \quad \left|\frac{e^{-i\kappa
r/\lambda}-e^{-i\kappa q/\mu}} {r/\lambda-q/\mu}\right|\le
e^{\varepsilon (r+q+m)}\|\kappa\|_{\spaceO},
$$
by Lemma \ref{lem:bound for compensator}.

The upper bounds  \eqref{eq:def f in J bounds} then imply that
\eqref{eq:acute a},\eqref{eq:acute b} converge in $\spaceO_{V\subset
\widetilde{U}}$, and
$$
\tilde{a}_{pq\l\mu}\le \tilde{C}\tilde{\rho}^{-p-q},\quad
\tilde{b}_{r\l}\le \tilde{C}\tilde{\rho}^{-r} \quad\text{ in
}\widetilde{U},
$$
which proves the first claim of the Lemma. The second claim is a
consequence of the first and of definition of $\M\spaceJ$-family. The
third claim follows from \eqref{eq:rotation} and \eqref{eq:shift}, so
the fourth follows from the third and the fact that $\operE_{\kappa}$
preserves $\M\spaceJ$.\end{proof}

\section{Petrov operator and proof of Theorem~\ref{thm:main}}\label{sec:petrov}
 We generalize the so-called
{\em Petrov operator}, see \cite{roitman-yakovenko}, as follows:
\begin{equation}\label{eq:Petrov oper}
\mathcal{P}_{\kappa}f(t):=\frac 1
{2i}\(f(te^{-i\kappa})-f(te^{i\kappa})\)\quad \kappa\in\R
\end{equation}

\begin{Lem}\label{lem:sin to Petrov}
The inverse Mellin transform of $g(s)\sin(\kappa s)$ is
$\mathcal{P}_{\kappa}f$. For any $f\in \spaceJ$
$$
\mathcal{P}_{\kappa}f=\M^{-1}\circ\left(\operE_{\kappa}
-\operE_{-\kappa}\right)\circ\M
f
$$
\end{Lem}
This is a straightforward applications of \eqref{eq:shift} and of
Lemma \ref{lem:e is well-defined} correspondingly.\qed

The key idea of the proof is a generalization of the Petrov trick
invented in \cite{petrov}, where it was used for $\kappa=\pi$.

Let $f(t)$  be a $\spaceJ$-function. Denote by $N(f)$ the number of
isolated zeroes of function $f$ on $(0,1)$ counted with
multiplicities, by $\Delta^1_{\kappa}(f)$ the increment of $\arg
f(t)$ along the arc $\{e^{i\kappa\phi}, \phi\in[-1,1]\}$ and let
$\Delta^0_\kappa(f)$ be equal to the limit as $\epsilon\to0$ of the
increment of $\arg f(t)$ along the arc $\{\epsilon e^{i\kappa\phi},
\phi\in[-1,1]\}$. Finiteness of all these numbers follows from
Lemma~\ref{lem:asympt for f}.

\begin{Lem}\label{lem:Petrov}
For $f\in \spaceJ$ and real on $(0,1)$ the following inequality
holds:
\begin{equation}\label{eq:Petrov}
N(f)\le 1+N(\mathcal{P}_{\kappa}f)+\frac 1
{2\pi}(\Delta^1_\kappa(f)+\Delta^0_\kappa(f)).
\end{equation}
\end{Lem}

\begin{proof}
Applying the argument principle to the sector $\{|\arg z|\le\kappa,
\epsilon\le|z|\le 1\}$ on the universal cover of the closed unit
disc we get
\begin{equation}
2\pi N(f)\le
\Delta_{-}+\Delta_{+}+\Delta_\kappa^1(f)+\Delta_\kappa^0(f),
\end{equation}
where $\Delta_{\pm}$ are  the increments of the argument of $f(t)$
along the rays $I_\pm=\{te^{\pm i\kappa}, t\in[\epsilon, 1]\}$.
Since $f(t)$ is real on the real axis, it takes conjugate values on
$I_\pm$, so $\Delta_{+}=\Delta_{-}$ and $\Im f|_{I_{\pm}}=\pm
\mathcal{P}_{\kappa}f$.

First, assume that $\mathcal{P}_{\kappa}f\not\equiv 0$. Evidently,
any segment of $I_+$ where increment of argument of $f(t)$ is
greater than $\pi$ contains at least one zero of $\Im f(t)$, so
$$\Delta_{+}\le \pi\left( N(\Im
f|_{I_+}+1\right)=\pi N(\mathcal{P}_{\kappa}f)+\pi,
$$
and we get the \eqref{eq:Petrov}.

If $\mathcal{P}_{\kappa}f\equiv 0$, then, by symmetry, $f(te^{
i\kappa})=f(te^{- i\kappa})$, i.e. $f=\tilde{f}(t^{\pi/\kappa})$
with $\tilde{f}$ holomorphic in a punctured unit disc. Therefore
\eqref{eq:Petrov} follows from the argument principle applied to
$\tilde f$ and a punctured unit disc.
\end{proof}

\begin{Rem} As $\kappa\to0$, $P_{\kappa}f(t)$ tends to a
multiple of $f'(t)$, the contributions $\Delta_\kappa^1(f)$ and
$\Delta_\kappa^0(f)$ of arcs tend to zero, and we get the standard
Rolle theorem as a limit case of Lemma~\ref{lem:Petrov}.
\end{Rem}

Let $\mathrm{F}=\{f(t,\nu)\}$ be a  $\spaceJ$-family on pair
$V\subset U$, and let $\kappa = \kappa(\nu)\in\spaceO_{V\subset U}$.

\begin{Lem}
The variation $\Delta_{\kappa(\nu)}^1 f$ is uniformly bounded from
above by some $\Delta^1=\Delta^1(\mathrm{F},\kappa)<\infty$ for all
$\nu\in V$.
\end{Lem}

{\em Proof.} The functions $\Im f(e^{i\kappa(\nu)\phi};\nu)$ and $\Re
f(e^{i\kappa(\nu)\phi};\nu)$ are analytic functions of
$\phi\in[-1,1]$ analytically depending on parameter $\nu\in
V\Subset\R^n$. Therefore their number of isolated zeros counted with
multiplicity is uniformly bounded from above by some number by
Gabrielov's theorem \cite{gabrielov}. The variation $\Delta^1$ is at
most $\pi(\text{this number}+1)$.\qed

\begin{Lem}
The variation $\Delta_{\kappa(\nu)}^0 f$ is uniformly bounded from
above by some $\Delta^0=\Delta^0(\mathrm{F},\kappa)<\infty$ for all
$\nu\in V$.
\end{Lem}

{\em Proof.} For a sum $f\in \spaceJ$  of a $\spaceJ$-series $\sigma$
the number $\Delta^0_\kappa(f)$ is bounded from above by
$-m_\sigma\kappa$ by Lemma \ref{lem:asympt for f}. Therefore the
existence of the uniform upper bound follows from the Lemma
\ref{lem:m is unif bdd in fam}.\qed

\begin{Cor}\label{cor:petrov-family}
For $\mathrm{F}$,  $\kappa$ as above and for any $f\in\mathrm{F}$
\begin{equation}\label{eq:Petrov-family}
N(f)\le N(P_{\kappa}f)+\Delta(\mathrm{F},\kappa),\quad
\Delta(\mathrm{F},\kappa)=\Delta^0(\mathrm{F},\kappa)+\Delta^1(\mathrm{F},\kappa)+1<\infty.
\end{equation}
\end{Cor}

\subsection{Proof of Theorem \ref{thm:main}} \label{ssec:proof}
Starting from the
analytic $\spaceJ$-family $\mathrm{F}$ of pseudoabelian integrals
parameterized by $(\theta,\omega)\in U\subset \mathscr{B}$ we
construct the new family $\{P_{\pi\l_n(\nu)}f(\bullet,\nu)\}$. It is
again an analytic $\spaceJ$-family by Lemma \ref{lem:e is
well-defined}. By Corollary \ref{cor:petrov-family}, application of
the operator $P_{\pi\l_n(\nu)}$ reduces the number of isolated zeros
at most by a finite number uniformly over parameters.

On the other hand, the new family is simpler: its Mellin transform
now has poles on $\{m/\l_k(\nu), k=1,...,n-1,
m\in\mathbb{Z}_{>-m_f}\}$, i.e. on the union of only $n-1$ arithmetic
progressions. Indeed, Mellin transform of
$\{P_{\pi\l_n(\nu)}f_{\nu}\}$ has the same poles as $\sin
(\pi\l_n(\nu) s)\M f$, and the first factor  has simple zeros at
$\{m/\l_n(\nu), m\in\mathbb{Z}_{>-m_f}\}$.

Repeating this step $n$ times we arrive at the  analytic
$\spaceJ$-family whose Mellin transform consists of entire functions.
This means that this analytic $\spaceJ$-family is identical zero by
Lemma \ref{lem:asympt for f}, so has no isolated zeros. \qed.

\end{document}